\documentclass[11pt]{article}
\usepackage{amsmath}
\usepackage{amssymb}
\usepackage{amsthm}
\usepackage[usenames]{color}
\usepackage{amscd}
\usepackage{dsfont}
\usepackage{indentfirst}
\usepackage[colorlinks=true,linkcolor=webgreen,filecolor=webbrown,
citecolor=webgreen]{hyperref}

\definecolor{webgreen}{rgb}{0,.5,0}
\definecolor{webbrown}{rgb}{.6,0,0}

\def\C{{\mathds{C}}}

\def\N{{\mathds{N}}}
\def\Z{{\mathds{Z}}}

\def\1{{\bf 1}}

\def\Res{\operatorname*{Res}}

\def\Hom{\operatorname{Hom}}

\newtheorem{theorem}{Theorem}[section]
\newtheorem{lemma}[theorem]{Lemma}
\newtheorem{cor}[theorem]{Corollary}

\newtheorem{prop}[theorem]{Proposition}

\begin{document}

\title{\bf Another generalization of the gcd-sum function}
\author{L\'aszl\'o T\'oth }
\date{}
\maketitle

\begin{abstract} We investigate an arithmetic function
representing a generalization of the gcd-sum function, considered by
Kurokawa and Ochiai in 2009 in connection with the multivariable
global Igusa zeta function for a finite cyclic group. We show that
the asymptotic properties of this function are closely connected to
the Piltz divisor function. A generalization of Menon's identity is
also considered.
\end{abstract}



\section{Introduction}

Let $r\in \N:=\{1,2,\ldots \}$ and define the arithmetic function
$A_r$ by
\begin{equation*} \label{def_A}
A_r(n):= \frac1{n^r} \sum_{k_1,\dots,k_r=1}^n \gcd(k_1\cdots k_r,n)
\quad (n\in \N).
\end{equation*}

The function $A_r$ was considered by Kurokawa and Ochiai
\cite{KurOch2009} in connection with certain zeta functions. More
exactly, the multivariable global Igusa zeta function for a group
$A$ is defined by
\begin{equation} \label{Igusa_def_gen}
Z^{\rm group}(s_1,\ldots,s_r;A):= \sum_{m_1,\ldots, m_r=1}^{\infty}
\frac{\# \Hom(A,\Z/m_1\cdots m_r\Z)}{m_1^{s_1}\cdots m_r^{s_r}}.
\end{equation}

Consider the case $A=\Z/n\Z$ ($n\in \N$). Since the number of group
homomorphisms $\Z/n\Z \to \Z/m_1\cdots m_r\Z$ is $\gcd(n,m_1\cdots
m_r)$, the function \eqref{Igusa_def_gen} reduces to
\begin{equation} \label{Igusa_def}
Z^{\rm group}(s_1,\ldots,s_r;\Z/n\Z):= \sum_{m_1,\ldots,
m_r=1}^{\infty} \frac{\gcd(m_1\cdots m_r,n)}{m_1^{s_1}\cdots
m_r^{s_r}}.
\end{equation}

Kurokawa and Ochiai \cite{KurOch2009} derived two representations
for \eqref{Igusa_def}, one of them being
\begin{equation} \label{Igusa_repr}
Z^{\rm group}(s_1,\ldots,s_r;\Z/n\Z)= \frac1{n^{s_1+\ldots+s_r}}
\sum_{k_1,\ldots, k_r=1}^{\infty} \gcd(k_1\cdots k_r,n)
\zeta(s_1,k_1/n) \cdots \zeta(s_r,k_r/n),
\end{equation}
where $\zeta(s,a):=\sum_{m=0}^{\infty} 1/(m+a)^s$ denotes the Hurwitz
zeta function. It follows from \eqref{Igusa_repr} that
\eqref{Igusa_def} has a meromorphic continuation to $\C^r$.

\begin{prop} \label{prop_formula_A_prime_power} {\rm (\cite[Cor.\ 1]{KurOch2009})}
For every $n=\prod_{p\mid n} p^{\nu_p(n)}\in \N$,
\begin{equation} \label{formula_A_prime_power}
A_r(n)= \prod_{p\mid n} \sum_{j=0}^r \left( \left( {\nu_p(n) \atop
j} \right) \right) \, \left(1-\frac1{p} \right)^j,
\end{equation}
where
\begin{equation*} \left(\left( {n \atop k} \right) \right):= \binom{n+k-1}{k}=
%\frac{n(n+1)\cdots (n+k-1)}{k!}
(-1)^k \binom{-n}{k}
\end{equation*}
denotes the number of $k$-multisets of an $n$-set.
\end{prop}

\begin{prop} {\rm (\cite[Cor.\ 2]{KurOch2009})} For every $n\in \N$,
\begin{equation} \label{limes_A} \lim_{r\to \infty} A_r(n)=n.
\end{equation}
\end{prop}

Formula \eqref{formula_A_prime_power} was obtained in
\cite{KurOch2009} as an application of the representations given for
\eqref{Igusa_def}, while \eqref{limes_A} is a direct consequence of
\eqref{formula_A_prime_power}. Note that
\eqref{formula_A_prime_power} was reproved in \cite{Min2009,Min2013}
using arguments of elementary probability theory.

In the case $r=1$,
\begin{equation} \label{def_A_1}
A_1(n):= \frac1{n} \sum_{k=1}^n \gcd(k,n)= \sum_{d\mid n}
\frac{\phi(d)}{d},
\end{equation}
where $\phi$ is Euler's totient function. Here $A_1(n)$ is
representing the arithmetic mean of $\gcd(1,n)$, $\ldots,\gcd(n,n)$
and \eqref{formula_A_prime_power} reduces to
\begin{equation*}
\label{formula_A_1_prime_power} A_1(n)= \prod_{p\mid n} \left( 1+
\nu_p(n) \left(1-\frac1{p} \right)\right).
\end{equation*}

See \cite{CheZha2011, KonKat2010, Tot2010, ZhaZha2011} for various
properties, analogs and other generalizations of the function
\eqref{def_A_1}.

In the present paper we derive a simple recursion formula for the
functions $A_r$, offer a direct number-theoretic proof for the
formula \eqref{formula_A_prime_power} and show that the asymptotic
properties of the function $A_r(n)$ are closely connected to the
Piltz divisor function $\tau_{r+1}(n)$, defined as the number of
ways of expressing $n$ as a product of $r+1$ factors.

As a modification of $A_r(n)$ we also consider and evaluate the
function
\begin{equation} \label{def_B_r}
B_r(n):= \sum_{\substack{k_1,\dots,k_r=1\\
\gcd(k_1\cdots k_r,n)=1}}^n \gcd(k_1\cdots k_r-1,n) \quad (n,r\in
\N).
\end{equation}

Note that in the case $r=1$,
\begin{equation} \label{def_B_1}
B_1(n):= \sum_{\substack{k=1\\
\gcd(k,n)=1}}^n \gcd(k-1,n) = \phi(n)\tau(n) \quad (n\in \N),
\end{equation}
where $\tau(n)$ stands for the number of divisors of $n$, according
to a result of Menon \cite{Kes1965}. See \cite{Tar2012,Tot2011M} for
other Menon-type identities.

Our results are given in Section \ref{Sect_2}, while their proofs
are included in Section \ref{Sect_3}.

\section{Results} \label{Sect_2}

Let $A_0(n):=\1(n)=1$ ($n\in \N$).

\begin{prop} \label{prop_recursion} The following recursion formula holds:
\begin{equation} \label{recursion_A}
A_r(n)= \sum_{d\mid n} \frac{\phi(d)A_{r-1}(d)}{d} \quad (n,r\in
\N).
\end{equation}
\end{prop}

Let $\overline{\phi}(n)=\phi(n)/n$.

\begin{cor} In terms of the Dirichlet convolution, $A_r
=\overline{\phi}A_{r-1}* \1$ ($r\in \N$). Therefore,
$A_1=\overline{\phi}* \1$, $A_2=\overline{\phi} (\overline{\phi}*
\1)* \1$, $A_3=\overline{\phi} (\overline{\phi} (\overline{\phi}
*\1)* \1)* \1$, in general
\begin{equation*}
A_r = \overline{\phi} (\overline{\phi} (\ldots (\overline{\phi}* \1)
\ldots )*\1)* \1
\end{equation*}
including $r$ times $\overline{\phi}$ and $r$ times $\1$.
\end{cor}

\begin{cor} \label{cor_multipl}
The function $A_r$ is multiplicative for any $r\in \N$.
\end{cor}

Observe that from formula \eqref{formula_A_prime_power},
\begin{equation} \label{ineq}
A_r(n)\le \prod_{p\mid n} \sum_{j=0}^r \binom{\nu_p(n)+j-1}{j}=
\prod_{p\mid n} \binom{\nu_p(n)+r}{r}= \tau_{r+1}(n)
\end{equation}
for any $n\in \N$, using parallel summation of the binomial
coefficients.

Also, $A_r(p^k) =\binom{k+r}{r}+{\cal O}(1/p) = \tau_{r+1}(p^k) +
{\cal O}(1/p)$, as $p\to \infty$ ($p$ prime) for any fixed $k,r\in
\N$. This suggests that the asymptotic behavior of $A_r(n)$ is
similar to that of $\tau_{r+1}(n)$.

\begin{prop} \label{prop_Dirichlet_series} The Dirichlet series of the function
$A_r$ has the representation
\begin{equation*}
\sum_{n=1}^{\infty} \frac{A_r(n)}{n^s}= \zeta^{r+1}(s) F_r(s) \quad
(\Re (s) > 1),
\end{equation*}
where the Dirichlet series $F_r(s):= \sum_{n=1}^{\infty} f_r(n)/n^s$
is absolutely convergent for $\Re (s) >0$. Moreover, for any prime
power $p^k$, $f_r(p^k)=0$ if $k\ge r+1$ and $f_r(p^k)\ll 1/p$, as
$p\to \infty$ if $1\le k\le r$.
\end{prop}

For the function $\tau_k$ ($k\ge 2$) one has
\begin{equation} \label{tau_k}
\sum_{n\le x} \tau_k(n) = \Res_{s=1}  x^s \frac{\zeta^k(s)}{s} +
\Delta_k(x),
\end{equation}
where the main term is $x P_{k-1}(\log x)$ with a suitable
polynomial $P_{k-1}(t)$ in $t$ of degree $k-1$ having the leading
coefficient $1/(k-1)!$. For the error term, $\Delta_k(x)= {\cal O}
(x^{\alpha_k+\varepsilon})$, with $\alpha_k\le (k-1)/(k+1)$ ($k\ge
2$), $\alpha_k\le (k-1)/(k+2)$ ($k\ge 4$). See \cite[Ch.\
XII]{Tit1986} and \cite{Ivi2006} for further results on
$\Delta_k(x)$.

\begin{prop} \label{prop_asymptotics} Let $r\in \N$. Then
\begin{equation} \label{asymp_A_r}
\sum_{n\le x} A_r(n) = x Q_r(\log x) + R_r(x),
\end{equation}
where $Q_r(t)$ is a polynomial in $t$ of degree $r$ having the
leading coefficient
\begin{equation*}
\frac1{r!} \prod_p \left( 1+ \sum_{k=1}^r
\frac{f_r(p^k)}{p^k}\right),
\end{equation*}
and $R_r(x)= {\cal O} (x^{\alpha_{r+1}+\varepsilon})$ (valid for
every $\varepsilon>0$).

Also, $R_r(x)=O(x^{r/(r+2)+\varepsilon})$ and $R_r(x)= \Omega
(b_r(x))$, where
\begin{equation*}
b_r(x)= (x\log x)^{\frac{r}{2r+2}} (\log_2 x)^{\frac{r+2}{2r+2}
((r+1)^{(2r+2)/(r+2)}-1)} (\log_3 x)^{-\frac{3r+2}{4r+4}},
\end{equation*}
$\log_j$ denoting the $j$-fold iterated logarithm.
\end{prop}

\begin{prop} \label{limsup_r} For every $r\in \N$,
\begin{equation} \label{limsup_A_r}
\limsup_{n\to \infty} \frac{\log A_r(n) \log \log n}{\log n}=
\log(r+1).
\end{equation}
\end{prop}

In the case $r=1$, formulae \eqref{asymp_A_r}, without the omega
result, and \eqref{limsup_A_r} were obtained by
Chi\-dam\-ba\-ras\-wa\-my and Sitaramachandrarao \cite[Th.\ 3.1,
4.1]{ChiSit1985}. In fact, both results were proved in
\cite{ChiSit1985} for a slightly more general function, namely for
$\psi_k(n)=\sum_{d\mid n} \phi_k(d)/d^k$, where $k\in \N$ and
$\phi_k(n)=n^k\prod_{p\mid n} (1-1/p^k)$ is the Jordan function of
order $k$. Here $A_1(n)=\psi_1(n)/n$.

For the function $B_r(n)$ defined by \eqref{def_B_r} we have
\begin{prop}  \label{menon} For every $n,r\in \N$,
\begin{equation*}
B_r(n)= \phi^r(n)\tau(n).
\end{equation*}
\end{prop}

\section{Proofs} \label{Sect_3}
%\begin{proof}

{\bf Proof of Proposition} \ref{prop_recursion}.
\begin{equation*}
A_r(n)= \frac1{n^r} \sum_{k_1,\dots,k_r=1}^n \sum_{d\mid
\gcd(k_1\cdots k_r,n)} \phi(d)= \frac1{n^r} \sum_{d\mid n} \phi(d)
\sum_{\substack{k_1,\dots,k_r=1 \\ k_1\cdots k_r\equiv 0 \text{ (mod
$d$)}}}^n 1,
\end{equation*}
where for fixed $k_1,\ldots,k_{r-1}$ the congruence $k_1\cdots
k_{r-1}k_r\equiv 0$ (mod $d$) has $\gcd(k_1\cdots k_{r-1},d)$
solutions $k_r$ (mod $d$) and has $(n/d)\gcd(k_1\cdots k_{r-1},d)$
solutions $k_r$ (mod $d$). Therefore,
\begin{equation} \label{1}
A_r(n)= \frac1{n^{r-1}} \sum_{d\mid n} \frac{\phi(d)}{d}
\sum_{k_1,\dots,k_{r-1}=1}^n \gcd(k_1\cdots k_{r-1},d),
\end{equation}
and writing $k_j=dq_j+s_j$ with $1\le s_j\le d$, $0\le q_j\le n/d-1$
($1\le j\le r-1$) we see that the inner sum is
\begin{equation*}
\sum_{\substack{1\le s_1,\ldots,s_{r-1}\le d \\ 0\le
q_1,\ldots,q_{r-1}\le n/d-1}} \gcd(s_1\cdots s_{r-1},d) =
\left(\frac{n}{d}\right)^{r-1} d^{r-1} A_{r-1}(d),
\end{equation*}
and inserting this into \eqref{1} we obtain \eqref{recursion_A}.
%\end{proof}
\newpage
\medskip
{\bf Proof of Proposition} \ref{prop_formula_A_prime_power}.

The function $n\mapsto A_r(n)$ is multiplicative by Corollary
\ref{cor_multipl}. Therefore, to obtain
\eqref{formula_A_prime_power} it is sufficient to consider the case
$n=p^k$ ($k\in \N$), a prime power. Let $x_r(k):=A_r(p^k)$ ($r\ge
0$) with a fixed prime $p$. From the recursion formula
\eqref{recursion_A} we have
\begin{equation*}
A_r(p^k) = 1+ \sum_{j=1}^k \left(1-\frac1{p}\right) A_{r-1}(p^j),
\end{equation*}
that is, by denoting $t:=1-1/p$,
\begin{equation} \label{x_r(k)}
x_r(k) = 1+ t \sum_{j=1}^k x_{r-1}(j)  \quad (r,k\in \N),
\end{equation}
where $x_0(k):=A_0(p^k) =1$ ($k\in \N$). Here $x_1(k) = 1+ t
\sum_{j=1}^k x_0(j)=1+kt$, $x_2(k) = 1+ t \sum_{j=1}^k x_1(j)= 1+ t
\sum_{j=1}^k (1+jt)= 1+kt+ \frac{k(k+1)}{2} t^2$, $x_3(k) = 1+ t
\sum_{j=1}^k x_2(j)= 1+ t \sum_{j=1}^k
\left(1+jt+\frac{j(j+1)}{2}t^2 \right) = 1+ kt + \frac{k(k+1)}{2}
t^2 + \frac{k(k+1)(k+2)}{6} t^3$.

We show by induction on $r$ that $x_r(k)$ is a polynomial in $t$ of
degree $r$ with integer coefficients which do not depend on $r$,
more exactly,
\begin{equation} \label{hyp_induction}
x_r(k) = 1+ \sum_{i=1}^r \left(\left( {k \atop i} \right) \right)
t^i.
\end{equation}

Assume that \eqref{hyp_induction} is valid for $r$. Then by
\eqref{x_r(k)} we obtain for $r+1$,
\begin{equation*}
x_{r+1}(k) = 1+ t \sum_{j=1}^k x_r(j) = 1+ t\sum_{j=1}^k \left(
1+\sum_{i=1}^r \left(\left( {j \atop i} \right) \right) t^i \right)
\end{equation*}
\begin{equation*}
= 1+kt+ \sum_{i=1}^r t^{i+1} \sum_{j=1}^k \binom{j+i-1}{i} = 1+
\sum_{i=0}^r \binom{k+i}{i+1} t^{i+1}
\end{equation*}
\begin{equation*}
= 1+ \sum_{i=1}^{r+1} \binom{k+i-1}{i} t^i = 1+ \sum_{i=1}^{r+1}
\left(\left( {k \atop i} \right) \right) t^i,
\end{equation*}
applying the upper summation formula. This completes the proof of
\eqref{formula_A_prime_power}.

\medskip

%\begin{proof}
{\bf Proof of Proposition} \ref{prop_Dirichlet_series}.

We use the conventions $\binom{a}{0}=1$ ($a\in \Z$),
$\binom{a}{b}=0$ ($a,b\in \N$, $a<b$). In terms of the Dirichlet
convolution, $A_r=\tau_{r+1}* f_r$, $f_r=A_r*\mu^{(r+1)}$ with
$\mu^{(r+1)}=\mu
* \cdots * \mu$ ($r+1$ times), where $\mu^{(r+1)}(p^k)= (-1)^k \binom{r+1}{k}$ for
any prime power $p^k$ ($k\in \N$).

Hence for any $k\in \N$,
\begin{equation*}
f_r(p^k)=\sum_{\ell=0}^k \mu^{(r+1)}(p^{\ell})A_r(p^{k-\ell})
%\end{equation*}
%\begin{equation*}
= \sum_{\ell=0}^k (-1)^{\ell} \binom{r+1}{\ell} \sum_{j=0}^r
\binom{j+k-\ell -1}{j}\left( 1-\frac1{p}\right)^j
\end{equation*}
\begin{equation} \label{f_r}
= \sum_{j=0}^r \left( 1-\frac1{p}\right)^j \sum_{\ell=0}^k
(-1)^{\ell} \binom{r+1}{\ell} \binom{j+k-\ell-1}{j},
\end{equation}
which is a polynomial in $1/p$ of degree $r$.

Here for any $k\ge r+1$,
\begin{equation*}
f_r(p^k)= \sum_{j=0}^r \left( 1-\frac1{p}\right)^j
\sum_{\ell=0}^{r+1} (-1)^{\ell} \binom{r+1}{\ell}
\binom{j+k-\ell-1}{j}=0,
\end{equation*}
since $\binom{j+k-\ell-1}{j}$ is a polynomial in $\ell$ of degree
$j$ and the inner sum is zero for any $0\le j\le r$ using the
identity
\begin{equation*} \label{power_ident}
\sum_{\ell=0}^{n} (-1)^{\ell} \ell^j \binom{n}{\ell}=0 \quad (0\le
j\le n-1).
\end{equation*}

Now for $1\le k\le r$ we obtain from \eqref{f_r} that the constant
term of the polynomial in $1/p$ giving $f_r(p^k)$ is
\begin{equation*}
c:= \sum_{j=0}^r \sum_{\ell=0}^k (-1)^{\ell} \binom{r+1}{\ell}
\binom{j+k-\ell-1}{j}
\end{equation*}
\begin{equation*}
= \sum_{\ell=0}^k (-1)^{\ell} \binom{r+1}{\ell} \sum_{j=0}^r
\binom{j+k-\ell-1}{j}
\end{equation*}
\begin{equation*}
= \sum_{\ell=0}^k (-1)^{\ell} \binom{r+1}{\ell} \binom{r+k-\ell}{r},
\end{equation*}
using parallel summation again.

Using now that $\binom{r+k-\ell}{r}= (-1)^{k-\ell}
\binom{-r-1}{k-\ell}$ we obtain
\begin{equation*}
c= (-1)^k \sum_{\ell=0}^k \binom{r+1}{\ell}
\binom{-(r+1)}{k-\ell}=0,
\end{equation*}
by Vandermonde's identity.

Therefore, $f_r(p^k)\ll 1/p$, as $p\to \infty$ for any $k\in
\{1,\ldots,r\}$. This shows that the Dirichlet series $F_r(s)$ is
absolutely convergent for $\Re (s) >0$.
%\end{proof}

\medskip
{\bf Proof of Proposition} \ref{prop_asymptotics}.
%\begin{proof}
Using Proposition \ref{prop_Dirichlet_series} and \eqref{tau_k} for
$k=r+1$,
\begin{equation*}
\sum_{n\le x} A_r(n) = \sum_{d\le x} f_r(d) \sum_{e\le x/d}
\tau_{r+1}(e)
\end{equation*}
\begin{equation*}
= \sum_{d\le x} f_r(d) \left(\frac{x}{d}P_r(\log (x/d)) +
\Delta_{r+1}(x/d) \right),
\end{equation*}
and \eqref{asymp_A_r} follows by usual estimates.

To obtain the omega result let $g_r$ denote the inverse under
Dirichlet convolution of the function $f_r$. Then $g_r$ is
multiplicative, $\tau_{r+1}=g_r*A_r$, so that
\begin{equation*}
\sum_{n\le x} \tau_{r+1}(n) = \sum_{d\le x} g_r(d) \sum_{e\le x/d}
A_r(e),
\end{equation*}
and the Dirichlet series $\sum_{n=1}^{\infty} g_r(n)/n^s$ is
absolutely convergent for $\Re (s)>0$. Now apply the $\Omega$-result
concerning the function $\tau_k$, due to Soundararajan
\cite{Sou2003}, for $k=r+1$. In the case $r=1$,
\begin{equation} \label{A_1_omega}
\sum_{n\le x} \tau(n)= \sum_{d\le x} \frac1{d} \sum_{e\le x/d}
A_1(e) = x\log x +(2\gamma- 1)x + \sum_{d\le x} \frac1{d} R_1(x/d) +
O(\log x).
\end{equation}

Assume that $R_1(x)=\Omega(b_1(x))$ does not hold. Then for every
$c>0$ there exists $x_c>0$ such that $|R_1(x)|\le c\, b_1(x)$ for
any $x\ge x_c$. Now inserting this into \eqref{A_1_omega}
contradicts that $\Delta(x)=\Omega(b(x))$. The same proof works out
also for $r\ge 2$.
%\end{proof}

\medskip

%\begin{proof}
{\bf Proof of Proposition} \ref{limsup_r}.

Similar to the proof of \cite[Th.\ 4.1]{ChiSit1985}. By
\eqref{ineq}, $A_r(n)\le \tau_{r+1}(n)$ ($n\in \N$). Therefore,
using that \eqref{limsup_A_r} holds for $\tau_{r+1}(n)$ instead of
$A_r(n)$ (\cite[Eq.\ 3.4]{SurSit1975}) we obtain that the given
$\limsup$ is $\le \log (r+1)$.

Furthermore, for squarefree $n$,
\begin{equation*}
A_r(n)=\prod_{p\mid n} \sum_{j=0}^r (1-1/p)^j = \prod_{p\mid n}
p\left(1-(1-1/p)^{r+1} \right)
\end{equation*}
\begin{equation*}
= \prod_{p\mid n} \left( r+1-\frac{r(r+1)}{2}\cdot
\frac1{p}+O(1/p^2)\right) = (r+1)^{\omega(n)} \prod_{p\mid n} \left(
1-\frac{r}{2}\cdot \frac1{p}+O(1/p^2)\right),
\end{equation*}
as $p\to \infty$ (for every fixed $r$).

Let $n_x=\prod_{x/\log x< p\le x} p$. Then
\begin{equation*}
\frac{\log A_r(n_x) \log \log n_x}{\log n_x}
\end{equation*}
\begin{equation*}
= \log (r+1) \frac{\omega(n_x)\log \log n_x}{\log n_x} + \frac{\log
\log n_x}{\log n_x} \log \prod_{p\mid n_x} \left( 1-\frac{r}{2}\cdot
\frac1{p}+O(1/p^2)\right).
\end{equation*}

By using familiar estimates, $\log n_x \sim x$, $\log \log n_x \sim
\log x$ and $\omega(n_x) \sim x/\log x$. Hence $\omega(n_x)\log \log
n_x/\log n_x \to 1$, as $x\to \infty$.

Also, $\prod_{p\le x} \left( 1-\frac{r}{2}\cdot
\frac1{p}+O(1/p^2)\right) \sim C_r/(\log x)^{r/2}$ with a suitable
constant $C_r$. Therefore, $\prod_{p\mid n_x} \left(
1-\frac{r}{2}\cdot \frac1{p}+O(1/p^2)\right) \to 1$ as $x\to
\infty$, and the result follows.

\medskip

%\begin{proof}
{\bf Proof of Proposition} \ref{menon}.

We use the following lemma, which follows easily by the
inclusion-exclusion principle, cf. \cite[Th.\ 5.32]{Apo1976}.

\begin{lemma} \label{lemma_phi} Let $n,d,x\in \N$ be such that $d\mid n$, $1\le x\le d$,
$\gcd(x,d)=1$. Then
\begin{equation*}
\# \{k\in \N: 1\le k \le n, k\equiv x \ \text{\rm (mod $d$)},
\gcd(k,n)=1 \}=\phi(n)/\phi(d).
\end{equation*}
\end{lemma}

We also need the following identity, which reduces to
\eqref{def_B_1} in the case $a=1$.

\begin{lemma} \label{lemma_Menon} Let $\gcd(a,n)=1$. Then
\begin{equation*}
\sum_{\substack{k=1\\ \gcd(k,n)=1}}^n \gcd(ak-1,n)=\phi(n)\tau(n)
\quad (n\in \N).
\end{equation*}
\end{lemma}

For the proof of Lemma \ref{lemma_Menon} write
\begin{equation*}
\sum_{\substack{k=1\\ \gcd(k,n)=1}}^n \gcd(ak-1,n)
=\sum_{\substack{k=1\\ \gcd(k,n)=1}}^n \sum_{d\mid \gcd(ak-1,n)}
\phi(d) = \sum_{d\mid n} \phi(d) \sum_{\substack{1\le k\le n\\
\gcd(k,n)=1\\ ak \equiv 1 \text{ (mod $d$)} }} 1,
\end{equation*}
and observe that for every $d\mid n$ the congruence $ak\equiv 1$
(mod $d$) has a unique solution (mod $d$), since $\gcd(a,n)=1$.
Therefore the inner sum is $\phi(n)/\phi(d)$ by Lemma
\ref{lemma_phi}. See also \cite[Cor.\ 14]{Tot2011M}.

Now for the proof of Proposition \ref{menon},
\begin{equation*}
B_r(n)= \sum_{\substack{k_1,\dots,k_{r-1}=1\\
\gcd(k_1\cdots k_{r-1},n)=1}}^n \sum_{\substack{k_r=1\\
\gcd(k_r,n)=1}}^n \gcd((k_1\cdots k_{r-1})k_r-1,n),
\end{equation*}
and applying Lemma \ref{lemma_Menon} for $a=k_1\cdots k_{r-1}$ we obtain that the
inner sum is $\phi(n)\tau(n)$. Hence,
\begin{equation*}
B_r(n)= \sum_{\substack{k_1,\dots,k_{r-1}=1\\
\gcd(k_1\cdots k_{r-1},n)=1}}^n \phi(n)\tau(n)= \phi^r(n)\tau(n).
\end{equation*}

\section{Acknowledgement} The author gratefully acknowledges support
from the Austrian Science Fund (FWF) under the projects Nr.
P20847-N18 and M1376-N18. The author thanks Professor Werner Georg
Nowak for very helpful discussions on the subject.

\medskip

L\'aszl\'o T\'oth  

Institute of Mathematics, Department of Integrative Biology 

Universit\"at f\"ur Bodenkultur, Gregor Mendel-Stra{\ss}e 33, A-1180
Wien, Austria 

and 

Department of Mathematics, University of P\'ecs 

Ifj\'us\'ag u. 6, H-7624 P\'ecs, Hungary 

E-mail: ltoth@gamma.ttk.pte.hu

\end{document}